\documentclass[a4paper,11pt]{article}
\usepackage[all]{xy}
\usepackage[english]{babel}
\usepackage[latin1]{inputenc}
\usepackage{amsfonts}
\usepackage{amsthm}
\usepackage{amsmath}
\usepackage{amsfonts}
\usepackage{latexsym}
\usepackage{amssymb}
\usepackage{mathrsfs}
\usepackage[usenames]{color}

\newtheorem{fed}{Definition}[section]
\newtheorem{teo}[fed]{Theorem}
\newtheorem{cor}[fed]{Corollary}
\newtheorem*{teo*}{Theorem}

\newtheorem{lem}[fed]{Lemma}
\newtheorem{pro}[fed]{Proposition}

\theoremstyle{definition}
\newtheorem{rem}[fed]{Remark}

\newtheorem{exa}[fed]{Example}

\newtheorem*{teoD}{Theorem (Douglas)}

\def\bdem{\begin{proof}}
\def\edem{\renewcommand{\qed}{\hfill $\blacksquare$}
\end{proof}}

\def\cH{\mathcal{H}}

\def\cH{\mathcal{H}}

\date{}
\begin{document}

\title{Proper splittings of Hermitian operators}

\author{Guillermina Fongi $^{a}$, M. Celeste Gonzalez $^{b}$ $^{c}$\\ 
\fontsize {9}{9} \selectfont{$^a$ 
Centro Internacional Franco Argentino de Ciencias de la Informaci\'on y de Sistemas, CIFASIS (CONICET-UNR)} \\ \selectfont \fontsize {9}{9} \selectfont{Ocampo y Esmeralda (2000)  Rosario, Argentina.}
\\
\fontsize {9}{9} \selectfont{$^b$ Instituto Argentino de Matem\'atica ``Alberto P. Calder\'on'', IAM-CONICET} \\ \selectfont \fontsize {9}{9} \selectfont{Saavedra 15, Piso 3 (1083), Buenos Aires, Argentina.}
\\
\fontsize {9}{9} \selectfont{$^c$ Instituto de Ciencias, Universidad Nacional de General Sarmiento, Argentina.}\\
\fontsize {9}{9} \selectfont{$^a$  gfongi@conicet.gov.ar, $^b$ celeste.gonzalez@conicet.gov.ar}}

\date{}
\maketitle


\selectfont {{\sl {Keywords: \ }}{Hermitian operators, positive orthogonal decompositions, splitting of operators, operator equations}

\begin{abstract}

In this article we deepen the study of proper splittings of Hilbert space operators, with special emphasis on proper splittings of Hermitian operators. On the one hand, 
 we improve characterizations  given in [Fongi $\&$ Gonzalez, J.
Math. Anal. Appl., 545 (2025) 129093] of the convergence of  both the polar proper and the Moore Penrose proper splittings. On the other hand, we introduce new proper splittings and we compare their convergence with those of the polar and the Moore Penrose  proper splittings.

\end{abstract}

\maketitle

\section{Introduction}

Splittings and proper splittings of matrices and operators defined on finite and infinite dimensional Hilbert spaces have long been studied  to obtain, by iterations, solutions of matrix equations or operator equations, respectively. The pioneer work on the treatment of splitting of matrices is \cite{varga}.   Numerous works followed that enriched the theory of the study of matrix and operator splittings, see \cite{MR348984,MR1628383,MR1113154,MR1286436,MR3671533, AriasGonzalezProperSplittings,FG-splitting} and references therein.

In this paper we deal with proper splittings of Hilbert space operators to approximate the Douglas' reduced solution of an operator equation. In the infinite dimensional Hilbert space setting, operator equations 
are often treated using tools such as variational formulations, regularization techniques for ill-posed problems, conjugate gradient  or  Krylov-type methods among others. The reader can resort to \cite{CarusoKrylov, Engl,  Gro, Magri} for references of this topics. 
 The framework of proper splittings yields an intrinsic operator theoretic methodology that is particularly well suited to closed range operators that may fail to be invertible.

A proper splitting of $T\in\mathcal{L}(\mathcal{H})$ (where $\mathcal L(\mathcal{H})$ is the algebra of bounded linear operators  defined on a Hilbert space $\mathcal H$) is a decomposition $T=U-V$, where $U,V\in\mathcal{L}(\mathcal{H})$ and $U$ has the same range and nullspace  as $T$. This kind of decomposition is applied to find the Douglas' reduced solution of a solvable operator equation $TX=S$, with $S\in\mathcal{L}(\mathcal{H})$, by means of the iterative process:
\begin{equation}\label{proceso iterativo introduccion}
X^{i+1}=U^\dagger VX^i +U^\dagger S,
\end{equation}
where $U^\dagger$ denotes the Moore Penrose inverse of $U$. 
 
It is known that the iteration (\ref{proceso iterativo introduccion}) converges if and only if $\rho(U^\dagger V)<1$ (where $\rho(\cdot)$ denotes the spectral radius). In this case (\ref{proceso iterativo introduccion}) converges to $T^\dagger S$, the Douglas' reduced solution of $TX=S$, see for example \cite{MR348984, AriasGonzalezProperSplittings, FG-splitting}. Reduced solutions of operator equations describe operators in many areas such as operator theory, sampling and reconstruction theory, approximation theory, electrical network theory, among others; see \cite{AndersonT1, AndersonT2, AriasCorachGonzalezGeneralizedInverses, Corach-MS, Fillmore-Williams, Massey Stojanoff,RDD-core}.

For every closed range operator $T\in\mathcal{L}(\mathcal{H})$ the polar proper splitting $T=U_T-V$ was defined in  \cite{AriasGonzalezProperSplittings} and \cite{FG-splitting}, where $U_T$ is the partial isometry of the polar decomposition of $T$. This proper splitting  converges if and only if $\|T\| <2$, see \cite{FG-splitting}.  This partition is advantageous  between other proper splittings because, when applying the iterative process (\ref{proceso iterativo introduccion}), the computation of the Moore-Penrose inverse of $U_T$ is reduced to taking the adjoint of $U_T$. In this article we deepen the study of this proper splitting.

The main purpose of this paper is to study proper splittings of Hermitian closed range operators. In this case, the positive orthogonal decomposition of a Hermitian operator will be a useful tool. For this class of operators, two particular proper splittings  were defined in \cite{FG-splitting}; namely the MP-proper splitting and the projection proper splitting. By means of the positive orthogonal decomposition we provide more detailed characterizations of these   splitting  and, in addition, we introduce new proper splittings.

The article is organized as follows. Section 2 introduces  notation and preliminary concepts that will be used along the work.  Section 3  is devoted to  the study of proper splittings of bounded linear operators. The main contributions here are Theorem \ref{UdaggerVpositivo}, which improves a sufficient condition for the convergence of general proper splittings previously established in \cite{FG-splitting}, and Proposition \ref{rho del polar}, which analyzes the rate of convergence of the polar proper splitting. 
In Section 4 we put  focus on proper splitting of Hermitian operators. The  central tool here is the orthogonal positive decomposition that a Hermitian operator admits. For this operator class we obtain characterizations of the convergence of  both the polar proper and the Moore-Penrose proper splittings by means of the positive orthogonal decomposition, Theorem \ref{maximo general}, Corollary \ref{maximo} and Theorem \ref{MP convergencia desacoplada}. Furthermore, we show that the projection proper splitting does not converge for nonpositive Hermitian operators. On the other hand, Theorem \ref{propersplittings desde la dpo} and Theorem \ref{proper desde dpo con Tdagger} establish new proper splittings derived from the orthogonal positive decomposition.
Finally, the article concludes with a study of proper splittings of the form
 $T=U-V$, with $U$ a positive operator.

\section{Preliminaries}

Throughout this article $\mathcal{H}$ denotes a complex Hilbert space with inner product $\langle \cdot, \cdot \rangle$ and $\mathcal{L(H)}$ is the algebra of bounded linear operators from $\mathcal{H}$ to $\mathcal{H}$. By $\|\cdot\|$ we denote the norm of an element in $\mathcal{H}$ induced by the inner product  $\langle \cdot, \cdot \rangle$ or the operator norm in $\mathcal{L}(\mathcal{H})$, according to the context. Given a densely defined operator $T$, $\mathcal D (T)$ denotes its domain. If $T\in\mathcal{L(H)}$ then  $\mathcal{R}(T)$ and $\mathcal{N}(T)$ stands for  the range and the nullspace of $T$, respectively. The adjoint operator of $T$ is denoted by $T^*$.  Recall that, $W\in\mathcal{L}(\mathcal{H})$ is a reflection if $W$ is an invertible operator such that  $W=W^{-1}$ and $W$ is a symmetry if $W=W^{-1}=W^{*} $.
In addition, $\mathcal{L}^h$ is the set of selfadjoint operators of $\mathcal{L(H)}$ and $\mathcal{L}^+$ is the set of positive operators of $\mathcal{L(H)}$.  Given $S,T\in\mathcal{L}^h$, we say that $S\leq T$ if $\langle Sx,x\rangle \leq \langle Tx,x\rangle$ for all $x\in\mathcal{H}$. The relation $\leq$ is the classic L\"owner operator order in $\mathcal{L}^h$.

Given $T\in \mathcal{L}(\mathcal H)$,   there exists a unique densely defined linear operator $T^\dagger$ with $\mathcal{D}(T^\dagger)=\mathcal{R}(T)\oplus\mathcal{R}(T)^\bot$  which solves, simultaneously, the four equations:
\begin{equation}\label{MP-equations}
TXT=T; \ \  XTX=X \textrm{ on } \mathcal R(T) \oplus \mathcal R(T)^{\bot}; \ \  XT = P_{T^*}\ \  TX=P_T|_{\mathcal R(T)\oplus \mathcal R(T)^\perp}.
\end{equation}
The operator $T^\dagger$ is called the Moore-Penrose inverse of $T$. It holds that $T^\dagger\in\mathcal{L}(\mathcal{H})$ if and only if $\mathcal{R}(T)$ is closed.

Given $T\in\mathcal{L}(\mathcal{H})$ we denote by $\sigma(T), \sigma_p(T), \sigma_{ap}(T)$ the spectrum, the point spectrum and the approximate point spectrum of $T$, respectively.  Remember that $\sigma_{ap}(T)\neq\emptyset$ for all $T\in\mathcal{L}(\mathcal{H})$ and $\sigma_p(T)\subseteq \sigma_{ap}(T)$.  Also it holds that  $\partial\sigma(T)\subseteq\sigma_{ap}(T)$ and if $T\in\mathcal{L}^h$ then  $\sigma_{ap}(T)\subseteq \sigma(T)$. By $\rho(T)$ we denote the spectral radius of $T$, i.e. $\rho(T)=\sup\{|\lambda|: \lambda\in\sigma(T)\}$. Recall that if $T \in \mathcal L(\mathcal H)$ is a normal operator, i.e., $TT^*=T^*T$, then  $\rho(T)=\|T\|$.

The following result, which characterizes the L\"owner order for positive operators, will be useful. Its proof is in   \cite[Proposition 2.4]{FG-splitting}.

\begin{lem}\label{BLT}
	Consider $S,T\in \mathcal L^+$ such that $\mathcal R(T)$ is closed. Then $S\leq T$ if and only if $\rho({T^\dagger S})\leq 1$ and $\mathcal{R}(S^{1/2})\subseteq \mathcal{R}(T)$.
\end{lem}
\medskip

The next result on range inclusion and factorization is due to Douglas \cite{MR0203464}:

\begin{teoD}
	Let $S,T \in \mathcal{L}(\mathcal{H})$. The following conditions are equivalent:
	\begin{enumerate}
		\item $\mathcal{R}(S)\subseteq \mathcal{R}(T)$;
		\item there exists a number $\lambda >0$ such that $SS^*\leq \lambda TT^*$;
		\item there exists $C\in\mathcal{L}(\mathcal{H})$ such that $TC=S$.
	\end{enumerate}
	In addition, if any of the above conditions holds then there exists a unique  $X_r\in\mathcal{L}(\mathcal{H})$ such that $TX_r=S$ and $\mathcal{R}(X_r)\subseteq \mathcal{N}(T)^\bot$. Furthermore, $\mathcal{N}(X_r)=\mathcal{N}(S)$ and  $X_r=T^ \dagger S$  is called the \it{Douglas' reduced solution} of $TX=S$.
\end{teoD} 

Given a closed subspace $\mathcal{S}\subseteq \mathcal{H}$, $P_{\mathcal{S}}$ denotes the orthogonal projection onto $\mathcal{S}$. Also, for $T\in \mathcal L(\mathcal H)$ we abbreviate $P_T$ to indicate the orthogonal projection onto $\overline{\mathcal{R}(T)}$. On the other hand, if  $\mathcal{S,T}\subseteq \mathcal{H}$ are two closed subspaces then the  orthogonal sum between $\mathcal{S}$ and $\mathcal{T}$ is denoted by  $\mathcal{S}\oplus\mathcal{T}$.

\smallskip

Recall that, given $T\in\mathcal{L(H)}$ there exists a unique partial isometry $U_T\in\mathcal{L(H)}$ with $\mathcal{N}(U_T)=\mathcal{N}(T)$ such that $T=U_T|T|=|T^*|U_T$, where $|T|=(T^*T)^{1/2}$. This factorization is called the polar decomposition of $T$.  
\smallskip

The reduced minimum modulus of a given $0\neq T\in\mathcal{L}(\cH)$ is defined by 
$$
\gamma({T})=\inf\{\|Tx\|: x\in\mathcal{N}(T)^\bot, \|x\|=1\}=\inf\sigma(|T|)\setminus\{0\}.
$$
It is well-known that $T$ has closed range if and only if $\gamma(T)>0$.  The following result will be useful along this article. Its proof is in \cite[Lemma 2.1]{Chium2021Onaconjeture}.

\begin{lem}\label{Lema Eduardo}
Let $T\in\mathcal{L}(\mathcal{H})$, $T\neq 0$ with polar decomposition $T=U_T|T|$. Then,
$$
\|T-U_T\|=\max\{1-\gamma(T), \|T\|-1\}.
$$
\end{lem}

\smallskip

We finish this  section with a collection of some results on the positive orthogonal decomposition of a Hermitian operator that will be useful throughout the  paper.
Given $T\in  \mathcal{L}^h$, the \textit{ positive orthogonal 
	decomposition} of $T$ is the unique  decomposition $T=T_1-T_2$ where $T_1, T_2\in \mathcal{L}^+$ and $\overline{\mathcal{R}(T_1)}\oplus\overline{\mathcal{R}(T_2)}=\overline{\mathcal{R}(T)}$.
It is well known that $T_1=\frac{|T|+T}{2}$ and $T_2=\frac{|T|-T}{2}$. Moreover, given $T \in \mathcal L^h$ and $T_1, T_2\in \mathcal L^+$, such that $T=T_1-T_2$   then  $T=T_1-T_2$ is the orthogonal decomposition of $T$ if and only if  $\mathcal{R}(T)=\mathcal{R}(T_1)\oplus\mathcal{R}(T_2)=\mathcal R(T_1+T_2)$. In particular,  $\mathcal{R}(T)$ is closed if and only if $\mathcal{R}(T_i)$ is closed, for $i=1,2$, see \cite[Lemma 4.2]{FM-desc.positivas} and \cite[Theorem 3.10]{AG-aditividad}. 

\medskip

\begin{rem}\label{dpo de MP}
   If $T=T_1-T_2$ is the positive orthogonal  decomposition of  a closed range operator $T\in \mathcal L^h$   then $T^\dagger =T_1^\dagger -T_2^\dagger.$ 
\end{rem}

\begin{pro} \label{dpo iso parcial}
Consider $U\in\mathcal{L}^h$. Then $U$ is a partial isometry if and only if $U=P_{\mathcal{S}}-P_{\mathcal{W}}$ is the positive orthogonal decomposition of $U$, where $\mathcal{S}$ and $\mathcal{W}$ are orthogonal closed subspaces. 
\end{pro}

\begin{proof}
Let $U\in\mathcal{L}^h$. If $U$ is a partial isometry then its positive orthogonal decomposition is $U=\frac{|U|+U}{2}-\frac{|U|-U}{2}=\frac{P_U+U}{2}-\frac{P_U-U}{2}$. Then to prove the assertion, it is sufficient to note that $\left(\frac{P_U+U}{2}\right)^2= \frac{P_U+U}{2}$ and  $\left(\frac{P_U-U}{2}\right)^2= \frac{P_U-U}{2}$. Conversely, if $U=P_{\mathcal{S}}-P_{\mathcal{W}}$, where $\mathcal{S}$ and $\mathcal{W}$ are ortohogonal  closed subspaces  then $U^2= P_{\mathcal{S}}+P_{\mathcal{W}}=P_{\mathcal{S}+\mathcal{W}}=P_U$, because $\mathcal R(U)=\mathcal S+\mathcal W$. Then $U$ is a partial isometry. 
\end{proof}

\begin{cor}\label{isoparcialautoadjunta}
Consider $T\in\mathcal{L}^h$ and  $T=T_1-T_2$ its positive orthogonal decomposition. If $T=U_T|T|$ is the polar decomposition of $T$ then
\begin{enumerate}
	\item $U_T=P_{T_1}-P_{T_2}$,
	\item $\|P_T-U_T\|=0$ or $\|P_T-U_T\|=2$.
\end{enumerate}  
\end{cor}

\begin{proof}
  1. By Proposition \ref{dpo iso parcial}, it holds that $P_{T_1}-P_{T_2}$ is a partial isometry. Observe that $\mathcal N(P_{T_1}-P_{T_2})=  \overline{\mathcal{R}(T_1)}\cap \overline{\mathcal{R}(T_2)}\oplus \mathcal N(T_1)\cap \mathcal N(T_2)=\mathcal N(T_1)\cap \mathcal N(T_2)=\mathcal N(T)$, where the first equality follows by \cite[Lemma 2.2]{KolRakFre2005}.  
 Finally, since $(P_{T_1}-P_{T_2})|T|=(P_{T_1}-P_{T_2})(T_1+T_2)=T_1-T_2=T$ then $U_T=P_{T_1}-P_{T_2}$. 
 
 2. Observe that if $T=T_1-T_2$ is the positive orthogonal decomposition of $T$, then $P_T-U_T=2P_{T_{2}}$. Therefore, $\|P_T-U_T\|=0$ or $\|P_T-U_T\|=2$.
\end{proof}

\section{Proper splittings of Hilbert space operators}

Given $T\in\mathcal{L}(\mathcal{H})$, a decomposition $T=U-V$, with $U,V\in\mathcal{L}(\mathcal{H})$,  is called a \textit{proper splitting} of $T$ if $\mathcal{R}(U)=\mathcal{R}(T)$ and $\mathcal{N}(U)=\mathcal{N}(T)$. This notion was introduced by Varga \cite{varga} in the matricial context in order to approximate the unique solution of the vectorial equation $Tx=s$ with $T$ an invertible matrix. Later, Berman and Plemmons \cite{MR348984}
 generalized this concept for rectangular matrices. In \cite{AriasGonzalezProperSplittings},  the notion of proper splitting of rectangular matrices  was considered in order to approximate reduced solutions of matrix equations. Then, in \cite{FG-splitting}, this analysis was extended to the context of Hilbert space operators. 
More precisely, proper splittings of $T$ are employed to obtain, by approximation,  the Douglas' reduced solution of a  solvable operator equation $TX=S$ through the iterative process
\begin{equation}\label{proceso iterativo}
X^{i+1}=U^\dagger VX^i +U^\dagger S.
\end{equation}

Of course we are interested in proper splittings $T=U-V$ with $V\neq 0$ because, if $V=0$,   there is no iteration  in process  (\ref{proceso iterativo}). 
It is known that the iterative process (\ref{proceso iterativo}) converges for all initial $X^0\in\mathcal{L(H)}$ if and only if $\rho(U^\dagger V)<1$.  Moreover, in this case, it converges to the  Douglas' reduced solution $T^\dagger S$ of $TX=S$, see \cite{AriasGonzalezProperSplittings, FG-splitting}.

\begin{rem}
The Douglas' reduced solution of operator equations appear in different context. Let us see some cases. Here, the symbol $\dot+$ denotes the direct sum between subspaces.
\begin{enumerate}
\item Many pseudinverses of a closed range operator can be defined as a reduced solution of an operator equation \cite{AriasCorachGonzalezGeneralizedInverses,RDD-core}. Namely, consider $T\in\mathcal{L}(\mathcal{H})$ with closed range:
\begin{enumerate}
\item  the Moore Penrose inverse of $T$ is the Douglas' reduced solution of the equation $TX=P_T$;
\item if $\mathcal{R}(T)\dot+\mathcal{N}(T)=\mathcal{H}$ then the group inverse of $T$ is the reduced solution for $\mathcal{R}(T)$ of the equation $TX=Q_{\mathcal{R}(T)//\mathcal{N}(T)}$;
\item if $\mathcal{R}(T)\dot+\mathcal{N}(T)=\mathcal{H}$ then the core inverse of $T$ is the reduced solution for $\mathcal{R}(T)$ of the equation $TX=P_T$;
\item if $\mathcal{R}(T)\dot+\mathcal{N}(T)=\mathcal{H}$ then the dual core inverse of $T$ is the Douglas' reduced solution of the equation $TX=Q_{\mathcal{R}(T)//\mathcal{N}(T)}$.
\end{enumerate}
It should be noted that for $T\in\mathcal{L}^h$ the group inverse, the core inverse and the dual core inverse coincides with the Moore Penrose inverse of $T$. It is necessary to mention that the iterative process \eqref{proceso iterativo} can be modified as process (2) in \cite{FG-splitting} in order to obtain a general reduced solution, see \cite{ AriasGonzalezProperSplittings, FG-splitting}  for details.

\item Douglas' reduced solutions are also relevant in reconstruction problems, in the context of frame theory. Let $\mathcal{U}$ be a closed subspace of $\mathcal{H}$. A sequence $\{u_n\}_{n\in\mathbb{N}}$ in $\mathcal{U}$ is called a frame for $\mathcal{U}$ if there exist numbers $A,B>0$ such that, for every $f\in\mathcal{U},$ it holds:
$$A||f||^2\leq \sum_{n\in\mathbb{N}}|<f,u_n>|^2\leq B||f||^2.$$
In order to describe a reconstruction problem, consider $f\in \mathcal{H}$ be the original (unknown) element and let $\{\langle f, u_j \rangle\}_{j\in\mathbb{N}}$ be a given sequence of samples, where $\{u_j\}_{j\in \mathbb{N}}$ is a frame for the subspace $\mathcal{U}=\overline{span}\{u_j\}_{j\in\mathbb{N}}$. 
The subspace $\mathcal{U}$ is called the {\it{sampling space}}. The synthesis operator of the frame $\{u_j\}_{j\in\mathbb{N}}$ is  $U\in\mathcal{L}(\ell^2,\mathcal{H})$, defined by $Uc=\sum_{n\in\mathbb{N}}c_n u_n.$ Then the sequence of samples is $U^*f$. On the other hand, let $\mathcal{T}$ be a closed subspace of $\mathcal{H}$, called the {\it{reconstruction space}}, and let $\{t_j\}_{j\in \mathbb{N}}$ be a frame of $\mathcal{T}$. The goal in a reconstruction problem is to approximate $f$ by a series expansion $\tilde{f}=\sum\limits c_j t_j\in \mathcal{T}$ obtained as $\tilde{f}=GU^*f$ for some $G\in \mathcal{L}(\ell^2,\mathcal{H})$ with $\mathcal{R}(G)\subseteq \mathcal{T}$.  The operator $G$ is called a {\it{reconstruction operator}}. In addition, $G$ is {\it{perfect}} if $GU^*g=g$ for all $g\in\mathcal{R}(T)$. It holds that there exists a perfect reconstruction operator $G$ if and only if $\mathcal{R}(T)\dot{+}\mathcal{N}(U^*)$ is direct and closed, where $T\in\mathcal{L}(\ell^2, \mathcal{H})$ is the synthesis operator of the reconstruction frame $\{t_j\}_{j\in\mathbb{N}}$. Under the hypothesis $\mathcal{R}(T) \dot{+}\mathcal{N}(U^*)$ is direct and closed, it holds that $G$ is a perfect reconstruction operator if and only if $GU^*$ is an idempotent operator with $\mathcal{R}(GU^*)=\mathcal{R}(T)$. Moreover, given an idempotent $Q\in\mathcal{L}(\mathcal{H})$ such that $\mathcal{R}(Q)=\mathcal{R}(T)$ and $\mathcal{N}(U^*)\subseteq \mathcal{N}(Q)$, the operator $G_Q=Q(U^*)^\dagger$ has optimal properties between all the possible perfect reconstruction operators $G$, see \cite[Section 3]{AriasGonzalezSampling} for details and proofs of these facts. Observe that $G_Q^*$ is the Douglas' reduced solution of the operator equation $UX=Q^*$.

\item  In some cases, the shorted operator of a Hermitian operator is also a Douglas' reduced solution of an operator equation. In fact, let $\mathcal{S}\subseteq \mathcal{H}$ a closed subspace and $A\in\mathcal{L}^h$ with closed range such that it is $P_{\mathcal S}$-positive (i.e., $(I-P_{\mathcal{S}})A(I-P_{\mathcal S})\in\mathcal{L}^+$). If, in addition, the pair $(A,\mathcal{S})$ is compatible (i.e., there exists some idempotent $Q\in\mathcal{L}(\mathcal{H})$ with $\mathcal{R}(Q)=\mathcal{S}$ which satisfies $AQ=Q^*A$) then the shorted operator $\Sigma(A,P_{\mathcal{S}})$ is the Douglas' reduced solution of the operator equation $A^\dagger X=P_AQ$ where $Q=I-E$, for $E$ an idempotent operator such that $AE=E^*A$ with $\mathcal R(E)=\mathcal{S}$. Also, it holds that $\mathcal{R}(\Sigma(A,P_{\mathcal{S}}))\subseteq \mathcal{R}(A)\cap\mathcal{S}$, see \cite{Massey Stojanoff} for details. 

\end{enumerate}
\end{rem}

\begin{lem}\label{caract de U}
	Consider $T\in \mathcal L(\mathcal H)$ and $T=U-V$, with  $U, V\in \mathcal L(\mathcal H)$,  $V\neq 0$. Then  $T=U-V$ is a  proper splitting of $T$ if and only if  $U=TY$ where $Y\in \mathcal L(\mathcal H)$ is such that  $\mathcal R(Y)=\overline{\mathcal R(T^ *)}	$, $ \mathcal N( Y)=\mathcal N(T)$ and $Y\neq P_{T^*}$.
		\end{lem}

\begin{proof}
Suppose that $T=U-V$ is a  proper splitting of $T$. Then the equation $TX=U$ is solvable and consider $Y$ the Douglas' reduced solution of $TX=U$. Therefore $\mathcal R(Y)\subseteq \overline{\mathcal R(T^*)}$ and $\mathcal N(Y) = \mathcal N (U)= \mathcal N(T)$. Observe that $\mathcal R(Y)=\mathcal R(T^ \dagger U)=T^ \dagger (\mathcal R(U))=T^ \dagger (\mathcal R(T))=\mathcal R(T ^ \dagger  T)=\overline{\mathcal R(T^*)}$. Conversely, if $U=TY$ where $Y\in \mathcal L(\mathcal H)$ is such that  $\mathcal R(Y)=\overline{\mathcal R(T^ *)}$ and $ \mathcal N( Y)=\mathcal N(T)$ then it is not difficult to see that $\mathcal R(U)=\mathcal R(T)$ and  $\mathcal N(U)=\mathcal N(T)$. So that $T=U-V$ is a proper splitting of $T$. Observe that $Y=P_{T^*}$ if and only if $U=T.$
\end{proof}

\begin{cor}
	Consider $T\in \mathcal L(\mathcal H)$ and $T=U-V$ a  proper splitting of $T$. Then $U^ \dagger V=P_{T^*}-Y^\dagger$, where $Y\in \mathcal L(\mathcal H)$ is the Douglas' reduced solution of $TX=U$.
\end{cor}

\smallskip

The following result provides a formula for the Moore-Penrose inverse of $T$ in terms of the elements of the proper splitting of $T$. This formula was given in \cite{MR348984} in the matricial context and in \cite{FG-splitting} in the closed range operators context.

\begin{lem}\label{Formula T dagger}
Let $T\in \mathcal{L}(\mathcal{H})$. If $T=U-V$ is a proper splitting of $T$ then $T^\dagger=(I-U^\dagger V)^{-1}U^\dagger$.
\end{lem}

\begin{proof}
Since $T=U-V$ is a proper splitting of $T$ then the equation $TX=U$ is solvable and, by Lemma \ref{caract de U} its reduced solution, $X_r=T^\dagger U$, satisfies that $\mathcal{R}(X_r)=\overline{\mathcal{R}(T^*)}$. Therefore, $X_r^\dagger=P_{T^*}X_r^\dagger=U^\dagger UX_r^\dagger=U^\dagger TX_rX_r^\dagger=U^\dagger TP_{T^*}=U^\dagger T$. Then $I-U^\dagger V=I-U^\dagger(U-T)=I-P_{T^*}+U^\dagger T=P_{ \mathcal{N}(T)}+X_r^\dagger$ is invertible. In fact, by \cite[Theorem 3.10]{AG-aditividad}, $\mathcal{R}(P_{ \mathcal{N}(T)}+X_r^\dagger)=\mathcal{H}$ and it is easy to check that $\mathcal{N}
(P_{ \mathcal{N}(T)}+X_r^\dagger)=\{0\}$. Note that $\mathcal{D}((I-U^\dagger V)^{-1}U^\dagger)=\mathcal{R}(U)\oplus\mathcal{R}(U)^\bot=\mathcal{D}(T^\dagger)$.  Now, we will check that $(I-U^\dagger V)^{-1}U^\dagger$ satisfies the equations (\ref{MP-equations}) in order to prove that $T^\dagger= (I-U^\dagger V)^{-1}U^\dagger$:
\begin{enumerate}
\item $T(I-U^\dagger V)^{-1}U^\dagger T= T(I-U^\dagger V)^{-1}(P_{T^*}-U^\dagger V)=T (P_{T^*}-U^\dagger V+P_{\mathcal{N}(T)})^{-1}(P_{T^*}-U^\dagger V)=T ((P_{T^*}-U^\dagger V)^\dagger +P_{\mathcal{N}(T)})(P_{T^*}-U^\dagger V)=T (P_{T^*}-U^\dagger V)^\dagger (P_{T^*}-U^\dagger V)=TP_{\mathcal{N}(P_{T^*}-U^\dagger V)^\bot}=TP_{T^*}=T,$
where the sixth equality follows because ${\mathcal{N}(P_{T^*}-U^\dagger V)}=\mathcal N(T).$ In fact, it is clear  that $\mathcal{N}(T)\subseteq \mathcal{N}(P_{T^*}-U^\dagger V)$. If $x\in \mathcal{N}(P_{T^*}-U^\dagger V)$ then $P_{T^*}x=U^\dagger Vx$. Then, by left multiplication by $U$ we get $Ux=Vx$, so that $Tx=(U-V)x=0$. Then $x\in\mathcal{N}(T)$.

\item$(I-U^\dagger V)^{-1}U^\dagger T(I-U^\dagger V)^{-1}U^\dagger=(I-U^\dagger V)^{-1}(P_{T^*}-U^\dagger V)(I-U^\dagger V)^{-1}U^\dagger=P_{T^*}(I-U^\dagger V)^{-1}U^\dagger=P_{T^*}((P_{T^*}-U^\dagger V)^\dagger+P_{\mathcal{N}(T)}) U^\dagger=(P_{T^*}-U^\dagger V)^\dagger U^\dagger= ((P_{T^*}-U^\dagger V)^\dagger+P_{\mathcal{N}(T)}) U^\dagger = (I-U^\dagger V)^{-1}U^\dagger$.

\item In item 1,  it has already been proven that $(I-U^\dagger 
V)^{-1}U^\dagger T=P_{T^*}$.

\item Finally, $\left(T(I-U^\dagger 
V)^{-1}U^\dagger\right)^2=T(I-U^\dagger V)^{-1}U^\dagger T(I-U^\dagger V)^{-1}U^\dagger=TP_{T^*}(I-U^\dagger V)^{-1}U^\dagger=T(I-U^\dagger V)^{-1}U^\dagger$. Furthermore, if $x\in\mathcal{N}(T)^\bot$ then $T(I-U^\dagger V)^{-1}U^\dagger Tx=TP_{T^*}x=Tx$. In addition, if $x\in\mathcal{R}(T)^\bot$ then $T(I-U^\dagger V)^{-1}U^\dagger x=0$. So that $T(I-U^\dagger V)^{-1}U^\dagger=P_{T}|_{\mathcal{R}(T)\oplus\mathcal{R}(T)^\bot}$.
\end{enumerate}
\end{proof}

The following sufficient condition for the convergence of general proper splittings improves \cite[Theorem 3.7]{FG-splitting} which was proved for certain closed range operators.

\begin{teo}\label{UdaggerVpositivo}
	Let $T\in\mathcal{L}(\mathcal{H})$ and $T=U-V$ a proper splitting of $T$.  If $T^\dagger V \in\mathcal{L}^+$ then $\rho(U^\dagger V)=\frac{\rho(T^\dagger V)}{1+\rho(T^\dagger V)}<1$. 
\end{teo}

\begin{proof}
	In the proof of Lemma \ref{Formula T dagger} we saw that if $X_r=T^\dagger U$ is the reduced solution of $TX=U$ then $\mathcal{R}(X_r)=\overline{R(T^*)}$ and $X_r^\dagger=U^\dagger T$. Now, suppose that $T^\dagger V \in\mathcal{L}^+$, then  $T^\dagger U \geq P_{T^*}$.  Therefore, by \cite[Theorem 4.1]{FG-MoorePenrose y orden} it holds that $U^\dagger T \leq P_{T^*}$, so that  $U^\dagger V \in\mathcal{L}^+$.
    
	If $T^\dagger V\in\mathcal{L}^+$ then $U^\dagger V=P_{T^*}-U^\dagger T\in\mathcal{L}^+$. In addition,  we assert that $\lambda\in \sigma_{ap}(T^\dagger V)$ if and only if 
	$\frac{\lambda}{1+\lambda}\in\sigma_{ap}(U^\dagger V)$. In fact, take
	$\lambda\in \sigma_{ap}(T^\dagger V)\subseteq \sigma(T^\dagger V)\subseteq[0,+\infty)$. Then there exists a sequence $\{x_n\}_{n\in\mathbb{N}}\subseteq \mathcal{H}$ such that $\|x_n\|=1$  and  $\|(T^\dagger V-\lambda I)x_n\|\underset{n\rightarrow \infty}{\longrightarrow} 0$. In addition, observe that $\|(T^\dagger V-\lambda I)x_n\|=\|((I-U^\dagger V)^{-1}U^\dagger V-\lambda I)x_n\|=\|(I-U^\dagger V)^{-1}(U^\dagger V-\lambda (I-U^\dagger V))x_n\|=\|(I-U^\dagger V)^{-1}(U^\dagger V(1+\lambda)-\lambda I)x_n\|=(1+\lambda)\|(I-U^\dagger V)^{-1}(U^\dagger V-\frac{\lambda}{1+\lambda} I)x_n\|\underset{n\rightarrow\infty}{\longrightarrow}0$. Since $(I-U^\dagger V)^{-1}(U^\dagger V-\frac{\lambda}{1+\lambda}I)x_n\underset{n\rightarrow\infty}{\longrightarrow}0$ if and only if $(U^\dagger V-\frac{\lambda}{1+\lambda}I)x_n\underset{n\rightarrow\infty}{\longrightarrow}0$, then the assertion follows. On the other hand, since $U^\dagger V\in\mathcal{L}^+$ then $\sigma_{ap}(U^\dagger V)\subseteq \sigma(U^\dagger V)\subseteq [0,+\infty)$. Then, since $\partial\sigma(U^\dagger V)\subseteq\sigma_{ap}(U^\dagger V)\subseteq \sigma(U^\dagger V)$, we get that
	\begin{eqnarray}
		\rho(U^\dagger V)&=&\sup\{\alpha:\alpha\in\sigma(U^\dagger V)\}=\sup\{\alpha: \alpha \in \partial\sigma(U^\dagger V)\}\nonumber\\
		&=& \sup\{\alpha: \alpha \in \sigma_{ap}(U^\dagger V)\}=\sup\left\{\frac{\lambda}{1+\lambda}: \lambda \in \sigma_{ap}(T^\dagger V)\right\}\nonumber\\
		&=&\frac{\rho(T^\dagger V)}{1+\rho(T^\dagger V)}<1. \nonumber
	\end{eqnarray}
\end{proof}

\begin{rem}
    In \cite[Theorem 3.7]{FG-splitting} the hypothesis $U^\dagger V\in \mathcal{L}^+\cap\mathcal{K}$ must be replaced by $T^\dagger V \in \mathcal{L}^+\cap\mathcal{K}$, because the positivity of $T^\dagger V$ implies the positivity of $U^\dagger V$ but the converse does not hold in general.  Note  that $T^\dagger V\in \mathcal{L}^+$ if $U^\dagger V\in \mathcal{L}^+$ and $\rho(U^\dagger V)<1$, see
 \cite[Proposition 6.1] {AriasGonzalezProperSplittings}.
    \end{rem}
\smallskip

The \textit{polar proper splitting} of a closed range operator $T\in\mathcal{L(\mathcal{H})}$ is $T=U_T-V$, where $U_T$ is the partial isometry of the polar decomposition of $T$, see \cite{AriasGonzalezProperSplittings} and \cite{{FG-splitting}}.
In \cite[Theorem 4.3]{FG-splitting} we proved the following equivalences for the convergence of this splitting. 

\begin{teo}\label{equivalencia polar con norma}
Let $T\in\mathcal{L}(\mathcal{H})$ with closed range and $T=U_T-V$ the polar proper splitting of $T$. The following assertions are equivalent:
\begin{enumerate}
\item The polar proper splitting of $T$ is convergent;
\item $\|T-U_T\|<1$;
\item $\|P_{T^*}-|T|\|<1$;
\item $\|T\|<2$.
\end{enumerate}
\end{teo}

We finish this section by giving more information about the polar proper splitting of a closed range operator.
If the operator equation $TX=S$ is solvable then the polar proper splitting of a multiple of $T$ can  be consider to obtain the reduced Douglas' solution of $TX=S$ using (\ref{proceso iterativo}), as we show in the following proposition.

\begin{pro}\label{alphaT para el polar}
Let $S,T\in\mathcal{L}(\mathcal{H})$ such that $T$ has closed range and the equation $TX=S$ is solvable. Also, let $\alpha\in\mathbb{R}^+$ such that $\|\alpha T\|<2$. Then, the polar proper splitting of $\alpha T$ is $U_T-W$ for some $W\in\mathcal{L}(\mathcal{H})$ and it converges to the reduced solution of the equation $TX=S$. 
\end{pro}

\begin{proof}
Let $T=U_T-V$ the polar proper splitting of $T$. Note that the polar decomposition of $\alpha T$ is $\alpha T=U_{\alpha T}|\alpha T|=U_{\alpha T}\;\alpha |T|$. Then,  $T=U_{\alpha T} |T|$ and so $U_{\alpha T}=U_T$. Note that the equations $\alpha TX=\alpha S$ and   $TX=S$  have the same the reduced solution. In addition, by Theorem \ref{equivalencia polar con norma}, the polar proper splitting of $\alpha T=U_T-W$ associated to the iterative method (\ref{proceso iterativo}) for the equation $\alpha TX=\alpha S$ converges. Then, the assertion follows.  \end{proof}

The following result provides a formula for  $\rho(U_T^*V)$ when the polar proper splitting of $T$ is convergent. This formula allows us to analyze the speed of convergence of the iterative process in terms of the norm of $T.$

\begin{pro}\label{rho del polar} Let $T\in\mathcal{L(\mathcal{H})}$ with closed range and $T=U_T-V$ the polar proper splitting of $T$. Then,
\begin{enumerate}
\item $\rho(U_T^*V)=1-\gamma(T)$,  if $0<\|T\|\leq1$;
\item $\rho(U_T^*V)=\|T\|-1$, if $1<\|T\|<2$.
\end{enumerate}
\end{pro}

\begin{proof}
Note that $0 < \gamma (T)\leq \|T||$. Also, by Lemma \ref{Lema Eduardo}, $\rho(U_T^*V)=\|U_T^*V\|=\|U_T-T\|=\max\{1-\gamma(T), \|T\|-1\}$. Then the assertions follow.
\end{proof}

\section{Proper splittings of Hermitian operators}

Hermitian operators have special properties that will be useful in the context of proper splittings.  For instance, if  $T\in\mathcal{L}^h$ has closed range then it is possible to get a proper splitting $T=U-V$ with $U\in\mathcal{L}^h$. Also, recall that  every Hermitian operator $T$ admits the canonical orthogonal decomposition $T=T_1-T_2$, where $T_1,T_2\in\mathcal{L}^+$ with $\mathcal{R}(T_1)$ and $\mathcal{R}(T_2)$ orthogonal. This decomposition simplifies the computation of $T^\dagger$, since   the sustractivity property holds, i.e. $T^\dagger=T_1^\dagger -T_2^\dagger$. In the following we will use the positive orthogonal decomposition of a Hermitian operator applied to $T$ or $U$, according to the convenience.

\subsection{The positive orthogonal decomposition on proper splittings.}

 In \cite[Section 6]{FG-splitting} we defined and studied the projection proper splitting and the Moore-Penrose proper splitting (MP-proper splitting) of a Hermitian operator. In what follows  we apply inherited techniques from the  positive orthogonal decompositions of Hermitian operators to deepen the study of proper splittings. In particular, this  decomposition allows us to give  more detailed descriptions of the convergences of the projection, the MP and the polar proper splittings. Moreover,  we introduce  new proper splittings of a Hermitian operator, by means of its positive orthogonal decomposition.

First, we show that even in the case that $T$ is an EP operator (i.e., $T$ with closed range and $\mathcal{R}(T)=\mathcal{R}(T^*)$), we can always consider $T=U-V$ a proper splitting  of $T$,  with $U\in\mathcal{L}^h$.

\begin{pro}
Let $T\in\mathcal{L}(\mathcal{H})$ be a closed range operator. Then, $T$ admits a proper splitting $T=U-V$ with $U\in\mathcal{L}^h$ if and only if $T$ is an EP operator.
\end{pro}

\begin{proof}
If $T$ is an EP operator it is sufficient to consider $U=P_T$. The converse is immediate.
\end{proof}

The next example illustrates a case where $T$ is an EP operator with $T=U-V$ a proper splitting  of $T$, where $U\in\mathcal{L}^h$ and   the computation of $U^\dagger$ is easier that the computation  of $T^\dagger$.
\begin{exa}
Consider $T\in\mathcal{L}(\ell^2)$ the EP operator defined by $T(x_1,x_2,x_3,\cdots)=(x_1+\frac{x_3}{2}, 0, \frac{x_1}{3}-x_3, x_4+\frac{x_6}{2},0,\frac{x_4}{3}-x_6, \cdots)$ and the proper splitting $T=U-V$ of $T$ where $U\in\mathcal{L}^h$ is  defined by $U(x_1,x_2,x_3,\cdots)=(x_1, 0, -x_3, x_4,0,-x_6, \cdots)$.
Then in this case it is advantageous to obtain the Douglas' reduced solution of a solvable $TX=S$   by approximation with the iterative process (\ref{proceso iterativo}) since $U^\dagger V$ is simpler than the calculation of $T^\dagger S$.
\end{exa}

\medskip

Now, let us  see that a solvable equation $TX=S$, where  $T\in\mathcal{L}^h$ is a closed range operator, can be uncoupled in two solvable equations by means of the positive orthogonal decomposition of $T$. Moreover, the Douglas' reduced solution of $TX=S$ is the sum of the Douglas' reduced solutions of the associated uncoupled equations.

\begin{pro} Let  $T\in\mathcal{L}^h$ with  closed range and $S\in\mathcal{L}(\mathcal{H})$ be such that $TX=S$ is a solvable equation.  If  $T=T_1-T_2$ is the positive orthogonal decomposition of $T$ and $X_r$ is the reduced solution of  $TX=S$ then  
$$
X_r=X_{r1}+X_{r2},
$$ where $X_{r1}$ is the Douglas reduced solution of  $T_1X=P_{T_1}S$ and $X_{r2}$ is the Douglas reduced solution of $T_2X=-P_{T_2}S$.
\end{pro}

\begin{proof}
Observe that $TX=S$ is solvable if and only if $T_1X=P_{T_1}S$ and $T_2X=-P_{T_2}S$ are solvable. In fact,  $\mathcal{R}(S)\subseteq \mathcal{R}(T)=\mathcal{R}(T_1)\oplus\mathcal{R}(T_2)$ if and only if $\mathcal{R}(P_{T_i}S)\subseteq \mathcal{R}(P_{T_i}T_i)=\mathcal{R}(T_i)$, $i=1,2$. Now, as $T^\dagger S=(T_1^\dagger-T_2^\dagger)S=T_1^\dagger P_{T_1}S-T_2^\dagger P_{T_2}S$. Then, the assertion follows. 
\end{proof}

Consider  a closed range operator $T\in \mathcal{L}^h$  and $T=U-V$ a proper splitting of $T$ with $U\in\mathcal{L}^h$. Then the iterative process (\ref{proceso iterativo})  can be uncoupled in terms of the positive orthogonal decomposition of $U$, $U=U_1-U_2$, as follows: 
\
\begin{eqnarray}
X^{i+1}& =& U^\dagger VX^i+U^\dagger S=(U_1^\dagger-U_2^\dagger)VX^i+(U_1^\dagger-U_2^\dagger)S\nonumber\\
&=&U_1^\dagger VX^i+U_1^\dagger S - U_2^\dagger VX^i+U_2^\dagger S.\label{proceso iterativo desacoplado}
\end{eqnarray}
If, in addition,  $T=T_1-T_2$ is the positive orthogonal decomposition of  $T$ and  $\mathcal{R}(U_j)=\mathcal{R}(T_j)$ then  each uncoupled iteration $U_j^\dagger VX^i+U_j^\dagger S$ comes from a proper splitting of  $T_j$; namely, $T_j=U_j-P_{U_j}V$, with $j=1,2$. Also, each iterative process converges  to the Douglas' reduced solution of the equation $T_jY=(-1)^{j+1}P_{T_j}S, $
with $j=1,2.$
Then, in this case,  it is natural to ask if there is any relationship between $\rho(U^\dagger V)$, $\rho(U_1^\dagger V)$ and $\rho(U_2^\dagger V)$.

\begin{teo}\label{maximo general}
Consider $T=T_1-T_2$ the positive orthogonal decomposition of a closed range operator $T\in \mathcal{L}^h$  and $T=U-V$ a proper splitting of $T$ with $U\in\mathcal{L}^h$ such that  the positive orthogonal decomposition of $U$, $U=U_1-U_2$, satisfies that $\mathcal{R}(U_i)=\mathcal{R}(T_i)$. If in addition  $U^\dagger T\in\mathcal{L}^h$, then the proper splitting  $T=U-V$ converges if and only if  $\rho(U^\dagger V)=\max\{\|U_1^\dagger V\|, \|U_2^\dagger V\|\}<1$.
\end{teo}

\begin{proof} 
Note that since $\mathcal{R}(U_i)=\mathcal{R}(T_i)$ then  $U^\dagger T\in\mathcal{L}^h$  if and only if  $U_i^\dagger T_i\in\mathcal{L}^h$, for $i=1,2$. In fact, $U^\dagger T= TU^\dagger $ if and only if $U_1^\dagger T_1 +U_2^\dagger T_2= T_1U_1^\dagger +T_2U_2^\dagger$; or equivalently $U_1^\dagger T_1 = T_1U_1^\dagger$ and $U_2^\dagger T_2=T_2U_2^\dagger$.
Observe that $\rho(U^\dagger V)=\|P_T-U^\dagger T\|=\|P_{T_1}+P_{T_2}-(U_1^\dagger -U_2^\dagger)T\|=\|P_{T_1}-U_1^\dagger T_1 + P_{T_2}-U_2^\dagger T_2 \|=\max\{\|(P_{T_1}-U_1^\dagger T_1\|,  \|P_{T_2}-U_2^\dagger T_2 \|\}$. Then the assertion follows. 
\end{proof}

The characterization of the convergence of the polar proper splitting of a Hermitian operator can be refined in terms   of the positive orthogonal decomposition of the operator.

\begin{cor}\label{maximo}
Consider a closed range operator $T\in \mathcal{L}^h$,  $T=T_1-T_2$ its  positive orthogonal decomposition and $T=U_T-V$ the polar proper splitting of $T$. Then the following assertions are equivalent:
\begin{enumerate}
\item The polar proper splitting of $T$ converges;
\item  $\max\{\|P_{T_1}-T_1\|, \|P_{T_2}-T_2\|\}<1$; 
\item $\|T_i\|<{2}$, for $i=1,2$.
\end{enumerate}
\end{cor}

\begin{proof}
$1. \leftrightarrow 2.$ Note that $U_T^ * T=|T|\in \mathcal L^h$. Then the equivalence  follows by Theorem \ref{maximo general}  and the fact that $U_T=P_{T_1}-P_{T_2} $.

 $1. \leftrightarrow 3.$ First observe that, since $\mathcal{R}(T_1)$ and $ {\mathcal{R}(T_2)}$ are orthogonal then  $\|T\|=\max\{ \|T_1\|, \|T_2\|\}$. Now,  by Theorem \ref{equivalencia polar con norma}, it holds that the polar proper splitting of $T$ converges if and only if $\|T\|<2$; or equivalently  $\|T_i\|<{2}$ for $i=1,2$. 
\end{proof}

\begin{rem}
Observe that the polar proper splitting is not the unique proper splitting which satisfies conditions of Theorem \ref{maximo general}.
In fact, consider $A\in\mathcal{L}^h$ with closed range, $A\neq0$ and $T\in\mathcal{L}^h(\mathcal{H}\oplus\mathcal{H}\oplus\mathcal{H}\oplus\mathcal{H})$ defined by 
$$T=\left(\begin{array}{cccc}
		\frac{1}{2} A&0& 0&0\\
		0&-\frac{2}{3} A&0 &0\\
        0&0& -\frac{1}{2}A & 0 \\
        0&0&0&0
        \end{array}\right), $$
and consider $T=U-V$  a proper splitting of $T$, where
$$U=\left(\begin{array}{cccc}
		A&0& 0&0\\
		0&-A &0 &0\\
        0&0& -3A & 0 \\
        0&0&0&0
        \end{array}\right) 
       \textrm{ and  } V=\left(\begin{array}{cccc}
		\frac{1}{2}A&0& 0&0\\
		0&-\frac{1}{3}A &0 &0\\
        0&0& -\frac{5}{2}A & 0 \\
        0&0&0&0
        \end{array}\right).$$

Let 
$T=T_1-T_2$
   and $U=U_1-U_2$
   be  the positive orthogonal decompositions of $T$ and  $U$, respectively. 
        Then it is not difficult to see that $U\neq U_T$, $U^\dagger T\in\mathcal{L}^h$ and also $\mathcal{R}(U_j)=\mathcal{R}(T_j)$, for $j=1,2$. 
	\end{rem}

Next, we  relate the convergence of the  polar proper splitting of a closed range Hermitian operator to the distance between the factors of the polar decomposition of the given operator.
\begin{pro}\label{convergencia del polar para positivos}
Let $T\in\mathcal{L}^h$ be a closed range operator and let $T=U_T|T|$ be its polar decomposition. Then $\||T|-U_T\|< 1$ if and only if $T\in \mathcal L^+$ and the polar proper splitting of $T$ converges.
\end{pro}

\begin{proof}
Suppose that $\||T|-U_T\|<1$. Since  $\||T|-U_T\|=\max\{\|T_1-P_{T_1}\|, \|T_2+P_{T_2}\|\}$ it follows that $T_2=0$, so that $T\in \mathcal L^+$. Hence,  by Lemma \ref{Lema Eduardo} it holds that  $\|T\|-1< 1$ and so, $\|T\|<2.$ Then, by Theorem \ref{equivalencia polar con norma} the polar proper splitting of $T$ converges. The converse follows also from Theorem \ref{equivalencia polar con norma}.
\end{proof}

For a normal closed range operator $T\in\mathcal{L}(\mathcal{H})$ the  \textit{projection proper splitting} of $T$ is $T=P_T-V$. It holds  that this splitting  converges if and only if $\|P_T-T\|<1$, see \cite{FG-splitting}. The technique of decomposing a Hermitian operator in its positive orthogonal decomposition allows us to see that this proper splitting  is not convergent for nonpositive Hermitian operators.

\begin{pro}\label{proj hermitiano no converge}
Consider $T\in\mathcal{L}^h\setminus\mathcal{L}^+$ a closed range operator. Then the projection proper splitting of $T$ does not converge. 
\end{pro}

\begin{proof}
Consider $T=T_1-T_2$ the positive orthogonal decomposition of $T$. Then $T-P_T=T_1-P_{T_1}-T_2-P_{T_2}$. Since $P_{T_2}\leq P_{T_2}+T_2$ then $\|T-P_T\|=\max\{\|T_1-P_{T_1}\|, \|T_2+P_{T_2}\|\}\geq 1$.
\end{proof}

 If $T\in \mathcal L(\mathcal{H})$ is a normal closed range operator then $T=T^\dagger -V$ is the \textit{MP-proper splitting} of $T$.  It holds that $T=T^\dagger -V$ converges if and only if $\|P_{T}-T^2\|<1$, see \cite{FG-splitting}. In the next result we state  a new characterization of the convergence of the MP-proper splitting of $T\in\mathcal{L}^h$ by means of its positive orthogonal decomposition.

\begin{teo}\label{MP convergencia desacoplada}
Consider $T\in \mathcal{L}^h$ with closed range and $T=T_1-T_2$ the positive orthogonal decomposition of $T$. Then, the following assertions are equivalent:
\begin{enumerate}
\item The MP-proper splitting of $T$ converges;
\item $\max\{\|P_{T_1}-T_1^2\|,\|P_{T_2}-T_2^2\|\}<1$;
\item $\|T_i\|<\sqrt{2}$, $i=1,2$;
\item $\|T\|<\sqrt{2}$.
\end{enumerate}
\end{teo}

\begin{proof}
Consider $T=T^\dagger-V$ the MP-proper splitting of $T$ and $T=T_1-T_2$ the positive orthogonal decomposition of $T$.

$1\leftrightarrow 2.$ Since $\rho(TV)=\rho(P_T-T^2)=\|P_T-T^2\|=\max\{\|P_{T_1}-T_1^2\|, \|P_{T_2}-T_2^2\|\}$ then the assertion follows.

$2\leftrightarrow 3.$ 
By Lemma \ref{Lema Eduardo}, for each $i=1,2$ it holds that $\|P_{T_i}-T_i^2\|=\max\{1-\gamma(T_i^2), \|T_i^2\|-1\}$. 
Since $\mathcal{R}(T_i)$ is closed then $1-\gamma(T_i^2)<1$. In addition, $\|T_i^2\|-1<1$ if and only if $\|T_i\|<\sqrt{2}$. Then the assertion follows. 

$3\leftrightarrow 4.$
Since $\|T\|=\max\{\|T_1\|, \|T_2\|\}$ the equivalence follows. 
\end{proof}

\begin{rem} Consider $S\in\mathcal{L}(\mathcal{H})$ and $T\in\mathcal{L}^h$ with closed range such that the equation $TX=S$ is solvable.
Similarly to the case of polar proper splittings (Proposition \ref{alphaT para el polar}),   the Douglas' reduced solution of the equation $TX=S$ can be obtained by means of a  MP-proper splitting without additional hypothesis on $T$. In fact, it suffices to consider the MP-proper splitting of  $\alpha T$, where  $\alpha\in\mathbb{R}^+$ is such that $\|\alpha T\|< \sqrt 2$.
\end{rem}

Given $T\in\mathcal{L}^h$,  there are proper splittings of $T$ that emerge from its positive orthogonal decomposition. However, not every convergent proper splitting is suitable for simplifying the computation of the Douglas' reduced  solution of a solvable equation $TX=S$. In the following example, we illustrate this assertion. 

\begin{exa}\label{splitting que no sirve}
Let $T\in\mathcal{L}^h$ and consider  $T=T_1-T_2$ the positive orthogonal decomposition of $T$.  Take $T=nT_1-mT_2 -V$ with $n,m\in \mathbb{N}$ and $n,m$ not simultaneously equal to 1. Note that $\mathcal{R}(nT_1-mT_2)=\mathcal{R}(T)$ and so that $\mathcal{N}(nT_1-mT_2)=\mathcal{N}(T)$. Then $T=nT_1-mT_2 -V$ is a proper division of $T$. Observe that
$T^\dagger V=(T_1^\dagger-T_2^\dagger) (nT_1-mT_2-T)=nP_{T_1}+mP_{T_2}
-P_T=(n-1)P_{T_1}+(m-1)P_{T_2}\in \mathcal L^+$. Hence, by Theorem \ref{UdaggerVpositivo}, it follows that the proper splitting $T=nT_1-mT_2 -V$ is convergent. Note that although this partition is convergent, the computation of $T^\dagger W$ by means of the iterative process $(\ref{proceso iterativo})$ has the same difficulty as the direct calculus of $T^\dagger W$. 
\end{exa}

Despite the previous example, there are other proper splittings of $T$, defined from its positive orthogonal decomposition, that could simplify the calculation of the Douglas' reduced solution employing the process $(\ref{proceso iterativo})$.

\begin{teo}\label{propersplittings desde la dpo}
    Consider a closed range operator $T\in \mathcal L^h$ with positive orthogonal decomposition $T=T_1-T_2$. Then the following assertions follow: \begin{enumerate}
    \item If $T_2\neq 0$ then $T=T_1-P_{T_2}-V$ is a proper splitting of $T.$ Moreover, $T=T_1-P_{T_2}-V$  converges if and only if $\|P_{T_2}-T_2\|<1.$
    \item If $T_1\neq 0$ then $T=T_2-P_{T_1}-V$ is a proper splitting of $T.$ Moreover, $T=T_2-P_{T_1}-V$  converges if and only if $\|P_{T_1}-T_1\|<1$.
    \end{enumerate}
\end{teo}

\begin{proof} $ \ $

1. First, note that $T=T_1-P_{T_2}-V$ is a proper splitting of $T.$ In fact,  it holds that  $\mathcal{R}(T_1-P_{T_2})=\mathcal{R}(T_1)\oplus\mathcal{R}(T_2)=\mathcal{R}(T)$. In addition, $\mathcal{N}(T_1-P_{T_2})=\mathcal{R}(T_1-P_{T_2})^\bot=(\mathcal{R}(T_1)\oplus\mathcal{R}(T_2))^\bot=\mathcal{R}(T)^\bot=\mathcal{N}(T)$.  
    Observe that $(T_1-P_{T_2})^\dagger V=(T_1-P_{T_2})^\dagger (T_1-P_{T_2}-T)=P_T-(T_1^\dagger- P_{T_2})(T_1-T_2)=P_T-P_{T_1}-T_2=P_{T_2}-T_2$. Hence, the last assertion holds.

    2. The proof is similar to the one  of  item 1.
    \end{proof}

In the next result, we compare the convergence of the polar proper splitting of a Hermitian operator with the splittings given in the above proposition.
\begin{cor}\label{maximo 2}
Consider a closed range operator $T\in \mathcal L^h$ with positive orthogonal decomposition $T=T_1-T_2$. The polar proper splitting $T=U_T-V$ of $T$ converges if and only if the proper splittings  $T=T_1-P_{T_2}-W$ and $T=T_2-P_{T_1}-Z$ converge. Moreover, in such case, $\max \{\rho((T_1-P_{T_2})^\dagger W), \rho((P_{T_1}-T_2)^\dagger Z)\}=\rho(U_T^*V)$. 
\end{cor}

\begin{proof}
The proof follows from  Theorem \ref{propersplittings desde la dpo} and Corollary \ref{maximo}.
\end{proof}

In the next we present a simple example, where we employ different splittings to approximate the shorted operator of a Hermitian operator  when applying the iterative method \eqref{proceso iterativo}. Moreover, in this case we  see  that one of the proper splittings of Theorem \ref{propersplittings desde la dpo} is more convenient than the polar proper splitting. 
\begin{exa}
Let $\mathcal{S}=\mathcal{S}_1\oplus\mathcal{S}_2\oplus\mathcal{S}_3\subseteq\mathcal{H}$ be a closed subspace, where $\mathcal{S}_i$ are closed subspaces for $i=1,\cdots, 3$. Let $A\in\mathcal{L}^h$ with closed range and an idempotent $Q\in\mathcal{L}(\mathcal{H})$ defined by 
$$
A=\left(\begin{array}{cccc}
2 I & 0 & 0 & 0\\
0& -3 I & 0  &0 \\
0 & 0 & 0 & 0 \\
0& 0& 0& 0
\end{array}\right)\begin{array}{c}
\mathcal{S}_1\\
\mathcal{S}_2\\
\mathcal{S}_3\\
\mathcal{S}^\bot
\end{array} \ \ \ \textrm{and } \ \ \ 
Q=\left(\begin{array}{cccc}
I & 0 & 0 & 0\\
0&  I & 0  & 0\\
0 & 0 & I & z \\
0& 0& 0& 0
\end{array}\right)\begin{array}{c}
\mathcal{S}_1\\
\mathcal{S}_2\\
\mathcal{S}_3\\
\mathcal{S}^\bot
\end{array},
$$
for some $z\in\mathcal{L}(\mathcal{S}^\bot,\mathcal{S}_3)$. Note that the pair $(A,\mathcal{S})$ is compatible because $AQ\in\mathcal{L}^h$. In order to approximate the shorted operator $\Sigma(A,P_\mathcal{S})$ we apply three different proper splittings of $A^\dagger$ to obtain the Douglas' reduced solution of $A^\dagger X=P_AQ$. Consider the positive orthogonal decomposition  $A^\dagger=A_1^\dagger-A_2^\dagger$ of $A^\dagger$. More precisely,
$$A^\dagger=\left(\begin{array}{cccc}
\frac{1}{2} I & 0 & 0 & 0\\
0& 0 & 0  &0 \\
0 & 0 & 0 & 0 \\
0& 0& 0& 0
\end{array}\right)- \left(\begin{array}{cccc}
0 & 0 & 0 & 0\\
0& \frac{1}{3}I & 0  &0 \\
0 & 0 & 0 & 0 \\
0& 0& 0& 0
\end{array}\right).$$
So, by Theorem \ref{propersplittings desde la dpo} we get that the proper splitting $A^\dagger=A_1^\dagger-P_{A_2^\dagger}-V$ converges and $\|P_{A_2^\dagger}-A_2^\dagger\|=\frac{2}{3}$. Also, the proper splitting $A^\dagger=A_2^\dagger-P_{A_1^\dagger}-W$ converges and $\|P_{A_1^\dagger}-A_1^\dagger\|=\frac{1}{2}$.
In addition, by Corollary \ref{maximo 2} we can see that the polar proper splitting of $A^\dagger$, $A^\dagger=U_{A^\dagger}-Z$ converges and $\rho(U_{A^\dagger}^*Z)=\frac{2}{3}$.
\end{exa}

\smallskip 

\begin{teo}\label{proper desde dpo con Tdagger}
    Consider a closed range operator $T\in \mathcal L^h$ with positive orthogonal decomposition $T=T_1-T_2$. Then the following assertions hold: \begin{enumerate}
    \item If $T_1\neq 0$ then $T=T_1^\dagger -T_2-V$ is a proper splitting of $T.$ Moreover, $T=T_1^\dagger -T_2-V$   converges if and only if the MP-proper splitting of $T_1$ converges.
      
    \item If $T_2\neq 0$ then $T=T_1 -T_2^\dagger-V$ is a proper splitting of $T.$ Moreover, $=T_1 -T_2^\dagger-V$   converges if and only if the MP-proper splitting of $T_2$ converges.
    \end{enumerate}
\end{teo}

\begin{proof} $ \ $

 1. Consider $U=T_1^\dagger -T_2$.   Observe that $\mathcal R(U)=\mathcal R( T_1^\dagger-T_2)=\mathcal R(T_1)\oplus \mathcal R(T_2)=\mathcal R(T).$ Similar to the proof of  Proposition \ref{propersplittings desde la dpo} item 1 it holds that  $\mathcal N(U)=\mathcal N(T).$ Finally, note that $U^\dagger V=(T_1- T_2^\dagger)(T_1^\dagger-T_2-T)=P_{T_1}-T_1^2$. Then the assertion follows.

 2. It is analogous to the above item.
\end{proof}

\begin{cor}\label{convergencia MP con nuevos}
Consider a closed range operator $T\in \mathcal L^h$ with positive orthogonal decomposition $T=T_1-T_2$. The MP proper splitting $T=T^ \dagger-V$ of $T$ converges if and only if the proper splittings  $T=T_1^ \dagger -{T_2}-W $and $T=T_1-{T_2}^\dagger-Z$ converge. Moreover, in such case, $\max \{\rho((T_1^ \dagger-{T_2})^\dagger W), \rho(({T_1}-T_2^\dagger)^ \dagger Z)\}=\rho(TV)$.
\end{cor}

\subsection{Splittings $T=U-V$, with $U\in\mathcal{L}^+$}

We have seen that if $T\in \mathcal L^h$ has closed range, then there exists a proper splitting of $T=U-V$ such that $U\in\mathcal L^h$. In fact, it is sufficient to consider $U=U_T.$ We ask now if for a  given $T\in \mathcal L^h$, it is possible to get a proper splitting with $U\in\mathcal{L}^+$.

\begin{pro}\label{splitting positivo}
Consider $T\in\mathcal{L}^h$ with closed range. Then $T=U-V$ is a proper splitting of $T$ if and only if  $T=UU^*-Z$ is a proper splitting of $T$ and $\mathcal{N}(T)=\mathcal{N}(U)$.
\end{pro}

\begin{proof}
Since $T=U-V$ is a proper splitting of $T$ then $\mathcal{R}(U)=\mathcal{R}(T)$ and $\mathcal{N}(U)=\mathcal{N}(T)$. Now, as $\mathcal{R}(UU^*)=\mathcal{R}(U)=\mathcal{R}(T)$ and $\mathcal{N}(UU^*)=\mathcal{N}(U^*)=\mathcal{R}(U)^\bot=\mathcal{R}(T)^\bot=\mathcal{N}(T)$ then we get that $T=UU^*-Z$ is a proper splitting of $T$. For the  converse it is sufficient to note that $\mathcal{N}(T)=\mathcal{N}(UU^*)=\mathcal{N}(U^*)$. Then $\overline{\mathcal R(U)}=\mathcal{R}(T)=\mathcal{R}(UU^*)\subseteq \mathcal{R}(U)$, so that $\mathcal{R}(T)=\mathcal{R}(U)$.
\end{proof}

Given $T\in \mathcal L^h$, 
it is natural to consider  the proper splitting $T=|T|-V$, because $\mathcal{R}(T)=\mathcal{R}(|T|)$. However, if $T=|T|-V$ converges then $\rho(|T|^\dagger V)=\|P_T-U_T\|<1$ so that $P_T=U_T$ and then  $T\in\mathcal{L}^+$, see  Corollary \ref{isoparcialautoadjunta}. This fact  was proved for closed range operators in \cite[Remark 6.3]{FG-splitting}.

\smallskip

Next, we  provide a formula to compute $\rho(U^\dagger V)$ for a proper splitting $T=U-V$ of $T\in\mathcal{L}^h$,  with $U\in\mathcal{L}^+$.

\begin{teo}\label{rho con norma}
Consider $T\in\mathcal{L}^h$ a  closed range operator and $T=U-V$  a proper splitting of $T$ such that $U\in \mathcal{L}^+ $. Then  $\rho(U^ \dagger V)=\|(U^{1/2})^\dagger (U-T)(U^{1/2})^ \dagger\|$. Moreover, if $\rho(U^ \dagger V)<1$ then $T\in \mathcal L^+$ and  $ U+V\in \mathcal L^+$.
\end{teo}

\begin{proof}   Observe that $\rho(U^\dagger V)=\rho((U^{1/2})^\dagger(U-T)(U^{1/2})^\dagger)=\|(U^{1/2})^\dagger(U-T)(U^{1/2})^\dagger\|$. Then the first assertion follows. 
Now, if $\rho(U^ \dagger V)<1$ then we get that  $\|(U^{1/2})^\dagger V(U^{1/2})^ \dagger\| <1$. Therefore  $-I<(U^{1/2})^\dagger V(U^{1/2})^ \dagger < I$, so that $-U\leq V\leq U$. Hence $0\leq T \leq 2U$ and so $T\in \mathcal L^+$ and $ U+V\in \mathcal L^+$.
\end{proof}

Note that the  above theorem generalizes \cite[Theorem 3.10]{MR3671533}.
Moreover, the converse of the second assertion is not true. In fact, it is sufficient to take 
$T=\left(\begin{array}{cc}
4&0\\
0&3\\
\end{array}\right) \in\mathbb{M}^{2\times 2}$ with $T=U-V$ a proper splitting of $T$,  where $U=\left(\begin{array}{cc}
1/2&0\\
0&1/2\\
\end{array}\right)$.

\begin{pro}\label{splitting con positivos}
Let $T\in\mathcal{L}^h$ with closed range and $T=U-V$ a proper splitting of $T$ with $U,V\in\mathcal{L}^+$. Then $\rho(U^\dagger V)\leq 1$ if and only if $T\in\mathcal{L}^+$.
\end{pro}

\begin{proof}
Let $T=U-V$ a proper splitting of $T$ with $U,V\in\mathcal{L}^+$. 
Note that $\mathcal{R}(V^{1/2})\subseteq \overline{\mathcal{R}(V)}=\overline{\mathcal{R}(T)}=\mathcal{R}(T)=\mathcal{R}(U)$. Then the assertion follows from Lemma \ref{BLT}.
\end{proof}

 Consider $T\in\mathcal{L}^h$ with closed range and $T=U-V=UU^*-Z$ proper splittings of $T$. In the following, we present some examples which illustrate that, in general,  there is  no relationship between the convergences of these splittings.

Theorem \ref{equivalencia polar con norma} and Proposition \ref{proj hermitiano no converge} show that given a closed range operator $T\in \mathcal L^h \setminus \mathcal L^+$ with $\| T\|<2$ then the polar proper $T=U_T-V$ is convergent but the projection proper splitting   
$T=U_TU_T^*-Z$ does not converge.
The next example provides another similar case.

\begin{exa}\label{ejemplo positivo 1}
Consider $T=T_1-T_2$ the positive orthogonal decomposition of a closed range operator $T\in\mathcal{L}^h\setminus \mathcal L^+$. We saw in Remark \ref{splitting que no sirve} that the proper splitting $T=nT_1-mT_2-V$ is always convergent. However, the proper splitting $T=(nT_1-mT_2)^2-W$ is not convergent.  In fact, since $((nT_1-mT_2)^2)^\dagger W=P_{T_1}+P_{T_2}-\frac{1}{n^2}T_1^\dagger+\frac{1}{m^2}T_2^\dagger$, so that  $\rho(((nT_1-mT_2)^2)^\dagger W)=\max{\{\|P_{T_1}-\frac{1}{n^2}T_1^\dagger\|,\|P_{T_2}+\frac{1}{m^2}T_2^\dagger\|\}}\geq 1$. Then the proper splitting $T=(nT_1-mT_2)^2-W$ is not convergent. 
\end{exa}

\begin{exa}\label{ejemplo positivo 2} Consider $A\in\mathcal{L}^h$ with closed range and $T\in\mathcal{L}^h(\mathcal{H}\oplus\mathcal{H}\oplus\mathcal{H})$ defined by 
$T=\left(\begin{array}{ccc}
I&0&0\\
0&A&0\\
0&0&0
\end{array}\right)$ and the proper splitting $T=U-V$, with $U=\left(\begin{array}{ccc}
2I&0&0\\
0&4A&0\\
0&0&0
\end{array}\right)$.
 Note that  $T=U-V$ is convergent because $\rho(U^\dagger V)<1$.
 Moreover,  the proper splitting $T=U^2-Z$ is also convergent because $\rho((U^2)^\dagger Z)<1$.
\end{exa} 

\begin{exa}\label{ejemplo positivo 3} Consider $A\in\mathcal{L}^h$ with closed range and $T\in\mathcal{L}^h(\mathcal{H}\oplus\mathcal{H}\oplus\mathcal{H})$ defined by 
$T=\left(\begin{array}{ccc}
I&0&0\\
0&3A&0\\
0&0&0
\end{array}\right)$. Then  the proper splitting $T=U-V$ with $U=\left(\begin{array}{ccc}
-\sqrt{2}I&0&0\\
0&-2A&0\\
0&0&0
\end{array}\right)$  is not convergent because $\rho(U^\dagger V)=\frac{5}{2}$.
However the proper splitting $T=UU^*-Z$ is convergent because $\rho((UU^*)^\dagger Z)<1$. 
\end{exa} 

\section{Concluding remarks}

In this work, we studied
proper splittings of bounded Hilbert space operators. 
For  operators in $\mathcal L(\mathcal H)$ we improved (Theorem \ref{UdaggerVpositivo}) a sufficient condition of convergence of  general proper splittings  given in  \cite[Theorem 3.7]{FG-splitting}.
For closed range operators, we gave a formula for the convergence's rate of the polar proper splitting in terms of the operator's norm and  another in terms of the operator's reduced minimum modulus, according to the norm of the operator (Proposition \ref{rho del polar}). 

On the other hand, we saw that the positive orthogonal decomposition of a Hermitian operator was a useful tool in the study of proper splittings of this class of operators. It allowed us to provide a new characterization of the convergence of the polar proper splitting and  the MP proper splitting  (Corollary \ref{maximo}, Theorem \ref{MP convergencia desacoplada}). In addition, we introduced two new proper splittings of an operator (Theorem \ref{propersplittings desde la dpo},  Theorem \ref{proper desde dpo con Tdagger}) originated from its positive orthogonal decomposition. We also compared their convergences with  those of the polar and MP proper splittings (Corollary \ref{maximo 2}, Corollary \ref{convergencia MP con nuevos}).  

Finally, we studied proper splittings $T=U-V$ of a Hermitian operator $T$ such that $U$ is a positive operator.
 Given $T=U-V$ a proper splitting of  a closed range Hermitian operator $T$,  it is always possible to consider the proper splitting $T = UU^* - W$ of $T$, where  $ UU^* $ is positive (Proposition \ref{splitting positivo}). However, we showed that there is no relationship between their convergences (Examples \ref{ejemplo positivo 1}, \ref{ejemplo positivo 2}, \ref{ejemplo positivo 3}).

\section*{Acknowledgements}
The authors thank the anonymous reviewers for their careful reading and their constructive comments that helped us improve the article.

\section*{Statements and Declarations}
 The authors have no conflicts of interest.

\end{document}